# On H-Irregularity Strength of Grid Graphs

Meilin Imelda Tilukay[1*]

[1] Department of Mathematics, Faculty of Mathematics and Natural Sciences, Pattimura University,
Jalan Ir. M. Putuhena, Kampus Poka – Unpatti, Ambon, Indonesia.
Email: meilin.tilukay@fmipa.unpatti.ac.id

**Abstract:** This paper deals with three graph characteristics related to graph covering named the (vertex, edge, and total, resp.) $H$–irregularity strength of a graph $G$ admitting $H$-covering. Those are the minimum values of positive integer $k$ such that $G$ has an $H$-irregular (vertex, edge, and total, resp.) $k$-labeling. The exact values of this three graph characteristics are determined for grid graph admitting grid-covering, .

## 1. Introduction

Let $G$ be a finite, simple, and undirected graph with the vertex set $V(G)$ and the edge set $E(G)$. A *labeling* of a graph is a mapping that sends some set of graph elements to a set of numbers (usually to positive of non-negative integer). If the domain is $V(G)$, or $E(G)$, or $V(G) \cup E(G)$, the labeling are called, respectively, a *vertex labeling*, or an *edge labeling*, or a *total labeling*. The most complete of recent survey on graph labeling given in [9].

The *weight of a vertex* $x \in V(G)$ under an edge $k$-labeling $\gamma: E(G) \to \{1, 2, \cdots, k\}$, where $k$ is a positive integer, is $w_\gamma(x) = \sum_{y \in E(G)} \gamma(xy)$, where the sum is over of all edges incident to $x$. The labeling $\gamma$ is called *irregular $k$-labeling* if all vertex-weights are pairwise distinct. The smallest positive integer $k$ for which $G$ has an irregular labeling is called the *irregularity strength* of $G$, denoted by $s(G)$. This graph characteristic was introduced by Chartrand, Jacobson, Lehel, Oellermann, Ruiz, and Saba in [7]. They gave the lower bound of $s(G)$ which is tight for many graphs given in [7, 8, 11]. Kalkowski, Karonski, and Pfender [11] also improved the bound of $s(G)$ to its best form.

The *weight of an edge* $xy \in E(G)$ under a vertex $k$-labeling $\delta: V(G) \to \{1, 2, \cdots, k\}$, where $k$ is a positive integer, is $w_\delta(xy) = \delta(x) + \delta(y)$. The labeling $\delta$ is called an *edge irregular $k$-labeling* if all edge-weights are pairwise distinct. The smallest positive integer $k$ for which $G$ has an edge irregular labeling is called the *edge irregularity strength* of $G$, denoted by $es(G)$. This graph characteristic was introduced by Ahmad, Al-Mushayt, and Baca in [1]. They gave the lower bound of $es(G)$ and proved that the lower bound is tight for some class of graphs.

Difficulties of finding the $s(G)$ or $es(G)$ for any graph or even for graphs with simple structure is one of motives how many researches on irregular labeling growing rapidly, including occurrence of its many variations, which can be seen in [ 9].

In [ 5], Baca, Jendrol, Miller, and Ryan combined the concepts of irregularity strength and a total labeling of a graph $G$ into a new characteristic of $G$, called the *total edge irregularity strength of* $G$, denoted by $tes(G)$. For a graph $G$, a total $k$-labeling $\varepsilon: V(G) \cup E(G) \to \{1, 2, \cdots, k\}$ is called an *edge irregular total $k$-labeling* of $G$ if for every pair of distinct edges $x_1 y_1$ and $x_2 y_2$, their weights are distinct, where the weight of an edge $xy$ under labeling $\varepsilon$ is $w_\varepsilon(xy) = \varepsilon(x) + \varepsilon(xy) + \varepsilon(y)$. The $tes(G)$ is the minimum $k$ for which $G$ has an edge irregular total labeling. They [5] provided the boundaries and proved its tightness for some certain graphs. Baca and Siddiqui [6] also gave the exact value of $tes(G)$ of generalized prism and Jendrol, Miskuf, and Sotak [10] for complete and complete bipartite graphs.

An *edge covering* of $G$ is a family of subgraphs $H_1, H_2, \cdots, H_t$ such that each edge of $E(G)$ belongs to at least one of the subgraphs $H_i, i = 1, 2, \cdots, t$. Then it is said that $G$ admits an $(H_1, H_2, \cdots, H_t)$-*(edge) covering*. If every subgraph $H_i$ is isomorphic to a given graph $H$, then $G$ admits an $H$-*covering*. Let $G$ be a graph admitting $H$-covering. In [2], Ashraf, Baca, Kimakova, and Semanicova-Fenovcikova introduced edge $H$-irregularity strength and vertex $H$-irregularity strength of a graph $G$, as a natural extensions of $s(G)$ and $es(G)$, with consideration of $H$-covering of $G$. For the subgraph $H \subseteq G$ under an edge $k$-labeling $\gamma: E(G) \to \{1, 2, \cdots, k\}$, the $H$-*weight* (i.e. the weight of a subgraph $H$) under labeling $\gamma$ is

$$w_\gamma(H) = \sum_{e \in E(H)} \gamma(e) \tag{1}$$

For the subgraph $H \subseteq G$ under a vertex $k$-labeling $\delta: V(G) \to \{1, 2, \cdots, k\}$, the $H$-*weight* under labeling $\delta$ is

$$w_\delta(H) = \sum_{v \in V(H)} \delta(v) \tag{2}$$

An edge $k$-labeling $\gamma$ is called an *$H$-irregular edge $k$-labeling* of a graph $G$ if for every two distinct subgraphs $H_1$ and $H_2$ isomorphic to $H$ there is $w_\gamma(H_1) \neq w_\gamma(H_2)$. The edge $H$-irregularity strength of a graph $G$, denoted by $ehs(G, H)$, is the smallest integer $k$ such that $G$ has an $H$ irregular edge $k$-labeling. In other hand, a vertex $k$-labeling $\delta$ is called an *$H$-irregular vertex $k$-labeling* of a graph $G$ if for every two distinct subgraphs $H_1$ and $H_2$ isomorphic to $H$ there is $w_\delta(H_1) \neq w_\delta(H_2)$. The vertex $H$-irregularity strength of a graph $G$, denoted by $vhs(G, H)$, is the smallest integer $k$ such that $G$ has an $H$ irregular vertex $k$-labeling. The lower bounds of $ehs(G)$ and $vhs(G)$ are given in [2] as follows.

*Theorem A.* Let $G$ be a graph admitting an $H$-covering given by $t$ subgraphs isomorphic to $H$. Then

$$2^{|E(G)|-1} \geq ehs(G, H) \geq \left\lceil 1 + \frac{t-1}{|E(H)|} \right\rceil.$$

*Theorem B.* Let $G$ be a graph admitting an $H$-covering given by $t$ subgraphs isomorphic to $H$. Then

$$2^{|V(G)|-1} \geq vhs(G, H) \geq \left\lceil 1 + \frac{t-1}{|V(H)|} \right\rceil.$$

The lower bound given in Theorem A is tight for path admits path and ladder admits cycle $C_4$-covering, while the lower bound given in Theorem B is tight for path admits path, ladder admits cycle, and fan admits cycle. They [2] also decrease the upper bounds of $ehs(G)$ and $vhs(G)$ in Theorem A and Theorem B for certain conditions.

Later, Ashraf, Baca, Lascsakova, and Semanicova-Fenovcikova in [3] introduced the total $H$-irregularity strength of a graph $G$, motivated by concepts the total irregularity strength, which is introduced by Marzuki, Salman, and Miller in [12] and $H$-covering of a graph $G$. A total $k$-labeling $f: V \cup E \to \{1, 2, \ldots, k\}$ of $G$ is called a totally irregular total $k$-labeling if for any pair of vertices $x$ and $y$, their weights $w(x)$ and $w(y)$ are distinct and for any pair of edges $x_1 x_2$ and $y_1 y_2$, their weights $w(x_1 x_2)$ and $w(y_1 y_2)$ are distinct. The minimum $k$ for which a graph $G$ has totally irregular total labeling, is called total irregularity strength of $G$, denoted by $ts(G)$. They [12] gave the lower bound of $ts(G)$ and provided the exact value of $ts(G)$ for cycles and paths. For fan, wheel, triangular book, and friendship graphs, Tilukay, Salman, and Persulessy [14] gave the exact values of their $ts$ which is equal to the lower bound, as well as Tilukay, Tomasouw, Rumlawang, and Salman gave in [15] for complete and complete bipartite graphs, while Ramdani and Salman [13] gave the exact values of $ts$ of some Cartesian products graphs.

Let $G$ be a graph admitting $H$-covering. For the subgraph $H \subseteq G$ under the total $k$-labeling $\varepsilon$, the $H$-weight under labeling $\varepsilon$ is defined as

$$w_\varepsilon(H) = \sum_{v \in V(H)} \varepsilon(v) + \sum_{e \in E(H)} \varepsilon(e) \qquad (3)$$

A total $k$-labeling $\varepsilon$ is called an *H-irregular total k-labeling* of $G$ if for every two distinct subgraphs $H_1$ and $H_2$ isomorphic to $H$, there is $w_\varepsilon(H_1) \neq w_\varepsilon(H_2)$. The *total H-irregularity strength* of $G$, denoted by $ths(G, H)$, is smallest integer $k$ such that $G$ has an $H$-irregular total $k$-labeling.

For the boundaries of $ths(G, H)$, they [3] gave the following result and proved the sharpness of the lower bound for path admits path, ladder admits cycle, and fan admits cycle $C_3$.

*Theorem C.* Let $G$ be a graph admitting an $H$-covering given by $t$ subgraphs isomorphic to $H$. Then

$$2^{|E(G)|-1} \geq ths(G, H) \geq \left\lceil 1 + \frac{t-1}{|V(H)| + |E(H)|} \right\rceil.$$

In [4], Ashraf, Baca, Semanicova-Fenovcikova, and Siddiqui proved the sharpness of the lower bound in Theorem A – C for ladders admit ladders and fans admit fans.

Note that for $H \cong K_2$, then for any graph $G$ admitting $K_2$-covering, $ths(G, K_2) = tes(G)$.

In this paper, the exact values of (vertex, edge, total) $H$-irregularity strength of grid graph $P_m \times P_n$. For $P_m \times P_n$ is Cartesian product of two paths $P_m$ and $P_n$, without losing generality, we set $m \leq n$.

## 2. Main Results

The first result gives the exact value of $vhs(P_m \times P_n, P_m \times P_c)$, with $2 \leq m \leq c \leq n$.

### 2.1. The Vertex Grid-Irregularity strength of Grid

*Theorem 1.* Let $P_m \times P_n$ and $P_m \times P_c$ be two grids, $2 \leq m \leq c \leq n$. Then

$$vhs(P_m \times P_n, P_m \times P_c) = \left\lceil 1 + \frac{n-c}{mc} \right\rceil.$$

*Proof.* Let $P_m \times P_n$ be a grid graph with the vertex set $V(P_m \times P_n) = \{u_i^j | 1 \leq i \leq m, 1 \leq j \leq n\}$ and the edge set $E(P_m \times P_n) = \{u_i^j u_{i+1}^j | 1 \leq i \leq m-1, 1 \leq j \leq n\} \cup \{u_i^j u_i^{j+1} | 1 \leq i \leq m, 1 \leq j \leq n-1\}$. It is a routine procedure to check that for every positive integer $c$, where $2 \leq m \leq c \leq n$, $P_m \times P_n$ admits $P_m \times P_c$ – covering with exactly $n - c + 1$ subgraphs isomorphic to $P_m \times P_c$. Since $|V(P_m \times P_c)| = mc$, by Theorem B, we have $vhs(P_m \times P_n, P_m \times P_c) \geq \left\lceil 1 + \frac{n-c}{mc} \right\rceil$.

Let $k = \left\lceil 1 + \frac{n-c}{mc} \right\rceil$. For the reverse inequality, we prove by constructing a $P_m \times P_c$ – irregular vertex $k$ – labeling $\lambda: V(P_m \times P_n) \to \{1, 2, \cdots, k\}$ as follows.

$\lambda(u_i^j) = \left\lceil 1 + \frac{j - ic}{mc} \right\rceil$, $1 \leq i \leq m, 1 \leq j \leq n$.

It can be checked that the maximum label is $(u_1^n) = \left\lceil 1 + \frac{n-c}{mc} \right\rceil = k$, which imply that $\lambda$ is a vertex $k$ – labeling.

Next, we evaluate the weights of all grid graphs $P_m \times P_c$ which cover $P_m \times P_n$. Let $H_l$, $1 \leq l \leq n - c + 1$, be a subgraph of $P_m \times P_n$ isomorphic to $P_m \times P_c$.

$V(H_l) = \{u_i^j | 1 \leq i \leq m, l \leq j \leq l + c - 1, \}$;

By equation (2), we have $\{w_\lambda(H_l) | 1 \leq l \leq n - c + 1\} = \{mc, mc + 1, mc + 2, \cdots, mc + n - c\}$ which made a consecutive sequence on difference 1, which means there is no two grids isomorphic to $P_m \times P_c$ of the same weight, and $\lambda$ is a $P_m \times P_c$ – irregular vertex $k$ – labeling. Thus, $vhs(P_m \times P_n, P_m \times P_c) = k$.    □

For our next results, we use the same definition on vertex (and edge) – set of grid $P_m \times P_n$ as in the proof of Theorem 1. The exact value of $ehs(P_m \times P_n, P_m \times P_c)$, with $2 \leq m \leq c \leq n$, is given below.

## 2.2. The Edge Grid-Irregularity strength of Grid

*Theorem 2.* Let $P_m \times P_n$ and $P_m \times P_c$ be two grids, $2 \leq m \leq c \leq n$. Then
$$ehs(P_m \times P_n, P_m \times P_c) = \left\lceil 1 + \frac{n-c}{2mc - m - c} \right\rceil.$$

*Proof.* Let $P_m \times P_n$ be a grid graph with the vertex set $V(P_m \times P_n) = \{u_i^j | 1 \leq i \leq m, 1 \leq j \leq n\}$ and the edge set $E(P_m \times P_n) = \{u_i^j u_{i+1}^j | 1 \leq i \leq m-1, 1 \leq j \leq n\} \cup \{u_i^j u_i^{j+1} | 1 \leq i \leq m, 1 \leq j \leq n-1\}$. Since $|E(P_m \times P_c)| = 2mc - m - c$, by Theorem A, we have $ehs(P_m \times P_n, P_m \times P_c) \geq \left\lceil 1 + \frac{n-c}{2mc-m-c} \right\rceil$.

Let $k = \left\lceil 1 + \frac{n-c}{2mc-m-c} \right\rceil$.

For the reverse inequality, we prove by constructing a $P_m \times P_c$ – irregular edge $k$ – labeling $\lambda: E(P_m \times P_n) \to \{1, 2, \cdots, k\}$ as follows.
$\lambda(u_i^j u_{i+1}^j) = \left\lceil 1 + \frac{j-ic}{2mc-m-c} \right\rceil, 1 \leq i \leq m-1, 1 \leq j \leq n;$
$\lambda(u_i^j u_i^{j+1}) = \left\lceil 1 + \frac{j-(m+i-1)c+i}{2mc-m-c} \right\rceil, 1 \leq i \leq m, 1 \leq j \leq n-1.$

It can be checked that the maximum label is $\lambda(u_1^n u_2^n) = \left\lceil 1 + \frac{n-c}{2mc-m-c} \right\rceil = k$, which imply that $\lambda$ is an edge $k$ – labeling.

Next, we evaluate the weights of all grid graphs $P_m \times P_c$ which cover $P_m \times P_n$. Let $H_l$, $1 \leq l \leq n-c+1$, be a subgraph of $P_m \times P_n$ isomorphic to $P_m \times P_c$.
$E(H_l) = \{u_i^j u_{i+1}^j | 1 \leq i \leq m-1, l \leq j \leq l+c-1\} \cup \{u_i^j u_i^{j+1} | 1 \leq i \leq m, l \leq j \leq l+c-2\};$
By equation (1), we have
$\{w_\lambda(H_l) | 1 \leq l \leq n-c+1\} = \{2mc-m-c, 2mc-m-c+1, 2mc-m-c+2, \cdots, 2mc-m+n-2c\}$
which made a consecutive sequence on difference 1, which means there is no two grids isomorphic to $P_m \times P_c$ of the same weight, and $\lambda$ is a $P_m \times P_c$ – irregular edge $k$ – labeling. Thus, $ehs(P_m \times P_n, P_m \times P_c) = k$.    □

Our last result gives the exact value of $ths(P_m \times P_n, P_m \times P_c)$, with $2 \leq m \leq c \leq n$.

## 2.3. The Total Grid-Irregularity strength of Grid

*Theorem 3.* Let $P_m \times P_n$ and $P_m \times P_c$ be two grids, $2 \leq m \leq c \leq n$. Then
$$ths(P_m \times P_n, P_m \times P_c) = \left\lceil 1 + \frac{n-c}{3mc - m - c} \right\rceil.$$

*Proof.* Let $P_m \times P_n$ be a grid graph with the vertex set $V(P_m \times P_n) = \{u_i^j | 1 \leq i \leq m, 1 \leq j \leq n\}$ and the edge set $E(P_m \times P_n) = \{u_i^j u_{i+1}^j | 1 \leq i \leq m-1, 1 \leq j \leq n\} \cup \{u_i^j u_i^{j+1} | 1 \leq i \leq m, 1 \leq j \leq n-1\}$. Since $|V(P_m \times P_c)| = mc$ and $|E(P_m \times P_c)| = 2mc - m - c$, by Theorem C, we have $ths(P_m \times P_n, P_m \times P_c) \geq \left\lceil 1 + \frac{n-c}{3mc-m-c} \right\rceil$.

Let $k = \left\lceil 1 + \frac{n-c}{3mc-m-c} \right\rceil$.

For the reverse inequality, we prove by constructing a $P_m \times P_c$ – irregular edge $k$ – labeling $\lambda: V(P_m \times P_n) \cup E(P_m \times P_n) \to \{1, 2, \cdots, k\}$ as follows.

$\lambda(u_i^j) = \left\lceil 1 + \frac{j-ic}{3mc-m-c} \right\rceil, 1 \leq i \leq m, 1 \leq j \leq n;$

$\lambda(u_i^j u_{i+1}^j) = \left\lceil 1 + \frac{j-(m+i)c}{3mc-m-c} \right\rceil, 1 \leq i \leq m-1, 1 \leq j \leq n;$

$\lambda(u_i^j u_i^{j+1}) = \left\lceil 1 + \frac{j-(2m+i-1)c+i}{2mc-m-c} \right\rceil, 1 \leq i \leq m, 1 \leq j \leq n-1.$

It can be checked that the maximum label is $\lambda(u_1^n) = \left\lceil 1 + \frac{n-c}{3mc-m-c} \right\rceil = k$, which imply that $\lambda$ is a total $k$ – labeling.

Next, we evaluate the weights of all grid graphs $P_m \times P_c$ which cover $P_m \times P_n$. Let $H_l$, $1 \leq l \leq n-c+1$, be a subgraph of $P_m \times P_n$ isomorphic to $P_m \times P_c$.

$V(H_l) = \{u_i^j | 1 \leq i \leq m, l \leq j \leq l+c-1, \};$

$E(H_l) = \{u_i^j u_{i+1}^j | 1 \leq i \leq m-1, l \leq j \leq l+c-1\} \cup \{u_i^j u_i^{j+1} | 1 \leq i \leq m, l \leq j \leq l+c-2\};$

By equation (3), we have

$\{w_\lambda(H_l) | 1 \leq l \leq n-c+1\} = \{3mc-m-c, 3mc-m-c+1, 3mc-m-c+2, \cdots, 3mc-m+n-2c\}$

which made a consecutive sequence on difference 1, which means there is no two grids isomorphic to $P_m \times P_c$ of the same weight, and $\lambda$ is a $P_m \times P_c$ – irregular total $k$ – labeling. Thus, $ths(P_m \times P_n, P_m \times P_c) = k$. □

## 3. Conclusion

For $2 \leq m \leq c \leq n$, grid $P_m \times P_n$ has a $P_m \times P_c$ – irregular vertex (edge, total) labeling, respectively, with a vertex (edge, total) $P_m \times P_c$ – irregularity strength, respectively, equals to its lower bound.

## References

[1] A. Ahmad, O.B.S. Al-Mushayt, M. Bača. (2014) On edge irregularity strength of graphs. *Appl. Math. Comput.* 243, 607–610.

[2] F. Ashraf, M. Bača, Z. Kimáková, A. Semaničová-Feňovčiková. (2016). On vertex and edge H-irregularity strengths of graphs, *Discrete Mathematics, Algorithms, and Applications* 8 (4), 1650070(1-13).

[3] F. Ashraf, M. Bača, M. Lasscsáková, A. Semaničová-Feňovčiková. (2017) On H-irregularity strength of graphs, *Discussiones Mathematicae Graph Theory,* 37 (4), 1067-1078.

[4] F. Ashraf, M. Bača, A. Semaničová-Feňovčiková, M.K. Siddiqui. On H-irregularity strength of ladders and fan graphs, *AKCE International Journal of Graphs and Combinatorics,* In Press.

[5] M. Bača, S. Jendrol̀, M. Miller, J. Ryan. (2007) On irregular total labellings, *Discrete Math.* 307, 1378-1388.

[6] M. Bača, M K. Siddiqui. (2014) Total edge irregularity strength of generalized prism. *Appl. Math. Comput.* 235, 168-173.

[7] G. Chartrand, M. S. Jacobson, J. Lehel, O.R. Oellermann, S. Ruiz, F. Saba. (1988) Irregular networks. *Congr. Numer.* 64, 187-192.

[8] R.J. Faudree, J. Lehel. (1987) Bound on the irregularity strength of regular graphs. *Combinatorica.* 52, 247–256.

[9] J A. Galian. (2019) A dynamic survey of graph labeling. *Electronic Journal of Combinatorics.* 22, 1-535 Retrieved Feb. 8th , 2020, from https://www. combinatorics.org/files/Surveys/ds6/ds6v22-2019. pdf

[10] S. Jendrol, J. Miskuf, R. Sotak. (2010) Total edge irregularity strength of complete graphs and complete bipartite graphs. *Discrete Mathematics.* 310, 400-407.

[11] M. Kalkowski, M. Karonski, F. Pfender. (2011) A new upper bound for the irregularity strength of graphs. *SIAM J. Discrete Math.* 25 (3), 1319–1321.

[12] C. C. Marzuki, A. N. M. Salman, M. Miller. (2013) On the total irregularity strengths of cycles and paths. *Far East Journal of Mathematical Sciences,* 82 (1), 1-21.

[13] R. Ramdani, A. N. M. Salman, (2013) On the total irregularity strengths of some cartesian products graphs. *AKCE Int. J. Graphs Comb.* 10(2), 199-209.

[14] M. I. Tilukay, A. N. M. Salman, E. R. Persulessy. (2015) On the total irregularity strength of fan, wheel, triangular book, and friendship graphs. *Procedia Computer Science*, 74, 124-131.

[15] M.I., Tilukay, B.P., Tomasouw, F.Y. Rumlawang, A.N.M. Salman, (2017) The total irregularity strength of complete graphs and complete bipartite graphs. *Far East Journal of Mathematical Sciences,* 102 (2), 317-327.